\documentclass{amsart}
\usepackage{amssymb}

\usepackage{amsmath}
\usepackage{amssymb}
\usepackage{amsthm}

\newtheorem{thm}{Theorem}[section]

\newtheorem{lem}[thm]{Lemma}

\theoremstyle{definition}
\newtheorem{defn}[thm]{Definition}
\theoremstyle{remark}

\numberwithin{equation}{section}

\begin{document}
\title{Automorphisms of finite Abelian groups}%

\author[C.J. Hillar]{Christopher J. Hillar}
\address{Department of Mathematics, Texas A\&M University, College Station, TX 77843}
\email{chillar@math.tamu.edu}

\author[D.L. Rhea]{Darren L. Rhea}
\address{Department of Mathematics, University of California, Berkeley, CA 94720}
\email{drhea@math.berkeley.edu}


\maketitle

\section{Introduction}

In introductory abstract algebra classes, one typically encounters
the classification of finite Abelian groups \cite{lang}:

\begin{thm}\label{classificationabelian}
Let $G$ be a finite Abelian group.  Then $G$ is isomorphic to a
product of groups of the form \[H_p = \mathbb Z/p^{e_1}\mathbb Z
\times \cdots \times \mathbb Z/p^{e_n}\mathbb Z,\] in which $p$ is
a prime number and $1 \leq e_1 \leq \cdots \leq e_n$ are positive
integers.
\end{thm}

Much less known, however, is that there is a description of
Aut$(G)$, the automorphism group of $G$.  The first compete
characterization that we are aware of is contained in a paper by
Ranum \cite{ranum} near the turn of the last century.  Beyond
this, however, there are few other expositions. Our goal is to
fill this gap, thereby providing a much needed accessible and 
modern treatment.

Our characterization of $\text{Aut}(G)$ is accomplished in three
main steps.  The first observation is that it is enough to work
with the simpler groups $H_p$.  This reduction is carried out by
appealing to a fact about product automorphisms for groups with
relatively prime numbers of elements (Lemma
\ref{autproductlemma}). Next, we use Theorem \ref{endquotientthm}
to describe the endomorphism ring of $H_p$ as a quotient of a
matrix subring of $\mathbb Z^{n \times n}$. And finally, the units
$\text{Aut}(H_p) \subset \text{End}(H_p)$ are identified from this
construction.

As a consequence of our investigation, we readily obtain an explicit
formula for the number of elements of $\text{Aut}(G)$ for any
finite Abelian group $G$ (see also \cite{Pan}).

\section{Product Automorphisms}

Let $G = H \times K$ be a product of groups $H$ and $K$, in which
the orders of $H$ and $K$ are relatively prime positive integers.
It is natural to ask how the automorphisms of $G$ are related to
those of $H$ and $K$.

\begin{lem}\label{autproductlemma}
Let $H$ and $K$ be finite groups with relatively prime orders.
Then \[  \text{\rm{Aut}}(H) \times \text{\rm{Aut}}(K) \cong
\text{\rm{Aut}}(H \times K).\]
\end{lem}

\begin{proof}
We exhibit a homomorphism $\phi: \text{\rm{Aut}}(H) \times
\text{\rm{Aut}}(K) \to \text{\rm{Aut}}(H \times K)$ as follows.
Let $\alpha \in \text{\rm{Aut}}(H)$ and $\beta \in
\text{\rm{Aut}}(K)$. Then, as is easily seen, an automorphism
$\phi(\alpha,\beta)$ of $H \times K$ is given by
\[\phi(\alpha,\beta)(h,k) = (\alpha(h),\beta(k)).\] Let
$\text{id}_H \in \text{Aut}(H)$ and $\text{id}_K \in
\text{Aut}(K)$ be the identity automorphisms of $H$ and $K$,
respectively.  To prove that $\phi$ is a homomorphism, notice that
$\phi(\text{id}_H,\text{id}_K) = \text{id}_{H \times K}$ and that
\[\phi(\alpha_1 \alpha_2, \beta_1 \beta_2)(h,k) = (\alpha_1
\alpha_2(h),\beta_1 \beta_2(k)) = \phi(\alpha_1 , \beta_1
)\phi(\alpha_2 , \beta_2 )(h,k),\] for all $\alpha_1, \alpha_2 \in
\text{Aut}(H)$, $\beta_1, \beta_2 \in \text{Aut}(K)$, and $h \in
H, k \in K$.

We next verify that $\phi$ is an isomorphism.  It is clear that
$\phi$ is injective; thus we are left with showing surjectivity.
Let $n = |H|$, $m = |K|$, and write $\pi_H$ and $\pi_K$ for the
standard projection homomorphisms $\pi_H: H \times K \to H$ and
$\pi_K: H \times K \to K$.  Fix $\omega \in \text{Aut}(H \times
K)$, and consider the homomorphism $\gamma: K \to H$ given by
$\gamma(k) = \pi_H(w(1_H,k))$, in which $1_H$ is the identity
element of $H$.  Notice that $\{k^n: k \in K\} \subseteq \ker
\gamma$ since \[1_H = \pi_H(w(1_H,k))^n = \pi_H(w(1_H,k)^n) =
\pi_H(w(1_H,k^n)) = \gamma(k^n).\]  Also, since $m$ and $n$ are
relatively prime,
the set $\{k^n: k \in K\}$ consists of $m$ elements.
Consequently, it follows that $\ker \gamma = K$ and $\gamma$ is
the trivial homomorphism. Similarly, $\delta: H \to K$ given by
$\delta(h) = \pi_K(w(h,1_K))$ is trivial.

Finally, define endomorphisms of $H$ and $K$ as follows:
\[\omega_H(h) = \pi_H(\omega(h,1_K)), \ \omega_K(k) =
\pi_K(\omega(1_H,k)).\]  From this construction and the above
arguments, we have \[\omega(h,k) = \omega(h,1_K) \cdot
\omega(1_H,k) = (\omega_H(h),\omega_K(k)) =
\phi(\omega_H,\omega_K)(h,k)\] for all $h \in H$ and $k \in K$. It
remains to prove that $\omega_H \in \text{Aut}(H)$ and $\omega_K
\in \text{Aut}(K)$, and for this it suffices that $\omega_H$ and
$\omega_K$ are injective (since both $H$ and $K$ are finite). To
this end, suppose that $\omega_H(h) = 1_H$ for some $h \in H$.
Then $w(h,1_K) = (w_H(h),w_K(1_K)) = (1_H,1_K)$, so $h = 1_H$
by injectivity of $w$. A similar argument shows that $\omega_K \in
\text{Aut}(K)$, and this completes the proof.
\end{proof}

Let $p$ be a prime number.  The order of $H_p = \mathbb
Z/p^{e_1}\mathbb Z \times \cdots \times \mathbb Z/p^{e_n}\mathbb
Z$ is easily seen to be $p^{e_1+\cdots + e_n}$.  As $G$ is
isomorphic to a finite product of $H_p$ over a distinct set of
primes $p$, Lemma \ref{autproductlemma} implies that Aut$(G)$ is
simply the product of Aut$(H_p)$ over the same set of primes.  We
will, therefore, devote our attention to computing Aut$(H_p)$ for
primes $p$ and integers $1 \leq e_1 \leq \cdots \leq e_n$.

\section{Endomorphisms of $H_p$}

In order to carry out our characterization, it will be necessary
to give a description of $E_p =$ End$(H_p)$, the endomorphism ring
of $H_p$. Elements of $E_p$ are group homomorphisms from $H_p$
into itself, with ring multiplication given by composition and
addition given naturally by $(A+B)(h) := A(h) + B(h)$ for $A,B \in
\text{End}(H_p)$ and $h \in H_p$.  These rings behave much like
matrix rings with some important differences that we discuss
below.

The cyclic group $C_{p^{e_i}} = \mathbb Z/p^{e_i}\mathbb Z$
corresponds to the additive group for arithmetic modulo $p^{e_i}$,
and we let $g_i$ denote the natural (additive) generator for
$C_{p^{e_i}}$.  Specifically, these elements $g_i$ can be viewed
as the classes \[\overline{1} = \{x \in \mathbb Z : x \equiv 1 \
(\text{mod } p^{e_{i}})\}\] of integers with remainder $1$ upon
division by $p^{e_i}$.

Under this representation, an element of $H_p$ is a vector
$(\overline{h}_1,\ldots,\overline{h}_n)^T$ in which each
$\overline{h}_{i} \in \mathbb Z/p^{e_i}\mathbb Z$ and $h_i \in
\mathbb Z$ is an integral representative. With these notions in
place, we define the following set of matrices.

\begin{defn}
\[R_p = \left\{ (a_{ij}) \in \mathbb Z^{n \times n} :  p^{e_i-e_j} \
| \ a_{ij} \text{ for all $i$ and $j$ satisfying $1 \leq j \leq i \leq n$}\right\}.\]
\end{defn}

As a simple example, take $n= 3$ with $e_1 = 1$, $e_2 = 2$, and $e_3 = 5$.
Then
\[R_p = \left\{
\left[\begin{array}{ccc}b_{11} & b_{12} & b_{13} \\b_{21}p & b_{22} & b_{23} \\b_{31}
 p^4 & b_{32}p^3 & b_{33}\end{array}\right] : b_{ij} \in \mathbb Z \right\}.\]
In general, it is clear that $R_p$ is closed under addition and
contains the $n \times n$ identity matrix $I$.  It turns out that
matrix multiplication also makes this set into a ring as the
following lemma demonstrates.

\begin{lem}
$R_p$ forms a ring under matrix multiplication.
\end{lem}

\begin{proof}
Let $A = (a_{ij}) \in R_p$.  The condition that $p^{e_i-e_j} \ | \
a_{ij}$ for all $i \geq j$ is equivalent to the existence of a decomposition \[ A =
PA'P^{-1},\] in which $A' \in \mathbb Z^{n \times n}$ and $P =
\text{diag}(p^{e_1},\ldots,p^{e_n})$ is diagonal. In particular,
if $A,B \in R_p$, then $AB = (PA'P^{-1})(PB'P^{-1}) = PA'B'P^{-1}
\in R_p$ as required.
\end{proof}

Let $\pi_i: \mathbb Z \to \mathbb Z/p^{e_i}\mathbb Z$ be the
standard quotient mapping $\pi_i(h) = \overline{h}$, and 
let $\pi: \mathbb Z^n \to H_p$ be the homomorphism given
by \[\pi(h_1,\ldots,h_n)^T = (\pi_1(h_1),\ldots,\pi_n(h_n))^T =
(\overline{h}_1,\ldots,\overline{h}_n)^T.\] We may now give a
description of $E_p$ as a quotient of the matrix ring $R_p$.
In words, the result says that an endomorphism of $H_p$ is multiplication
by a matrix $A \in R_p$ on a vector of integer 
representatives, followed by an application of $\pi$.

\begin{thm}\label{endquotientthm}
The map $\psi: R_p \to \text{\rm{End}}(H_p)$ given by
\[\psi(A)(\overline{h}_1,\ldots,\overline{h}_n)^T =
\pi(A(h_1,\ldots,h_n)^T)\] is a surjective ring homomorphism.
\end{thm}
\begin{proof}
Let us first verify that $\psi(A)$ is a well-defined map from
$H_p$ to itself.  Let $A = (a_{ij}) \in R_p$, and suppose that
$(\overline{r}_1,\ldots,\overline{r}_n)^T =
(\overline{s}_1,\ldots,\overline{s}_n)^T$ for integers $r_i,s_i$
(so that $p^{e_i} \ | \ r_i-s_i$ for all $i$).  The $k$th vector
entry of the difference $\pi(A(r_1,\ldots,r_n)^T) -
\pi(A(s_1,\ldots,s_n)^T)$ is
\begin{equation}
\begin{split}
\pi_k \left( \sum\limits_{i = 1}^n {a_{ki} r_i } \right)  - \pi_k
\left( \sum\limits_{i = 1}^n {a_{ki} s_i } \right) = \ & \pi_k
\left( \sum\limits_{i = 1}^n {a_{ki} r_i }  -
\sum\limits_{i = 1}^n {a_{ki} s_i } \right) \\
= \ &  \sum\limits_{i = 1}^n { \pi_k  \left( \frac{{a_{ki} }}
{{p^{e_k - e_i } }} \cdot p^{e_k  - e_i } ( {r_i  - s_i }) \right) }  \\
= \ & \overline{0}, \\
\end{split}
\end{equation}
since $p^{e_k} \ | \ p^{e_k  - e_i }(r_i  - s_i )$ for $k \geq i$
and $p^{e_k} \ | \ (r_i  - s_i )$ when $k < i$.  Next, since
$\pi$ and $A$ are both linear, it follows that $\psi(A)$ is
linear.  Thus, $\psi(A) \in \text{End}(H_p)$ for all $A \in R_p$.


To prove surjectivity of the map $\psi$, let $w_i =
(0,\ldots,g_i,\ldots,0)^T$ be the vector with $g_i$ in the $i$th
component and zeroes everywhere else. An endomorphism $M \in
\text{End}(H_p)$ is determined by where it sends each $w_i$;
however, there isn't complete freedom in the mapping of these
elements. Specifically, suppose that $M(w_j) =
(\overline{h}_{1j},\ldots,\overline{h}_{nj})^T =
\pi(h_{1j},\ldots,h_{nj})^T$ for integers $h_{ij}$.  Then, \[ 0 =
M\left( 0 \right) = M\left( {p^{e_j } w_j } \right) = \underbrace
{Mw_j  + \cdots + Mw_j }_{p^{e_j } } = \left(\overline{p^{e_j
}h_{1j}},\ldots,\overline{p^{e_j }h_{nj}} \right)^T.\]
Consequently, it follows that $p^{e_i} \ | \ p^{e_j}h_{ij}$ for
all $i$ and $j$, and therefore $p^{e_i-e_j} \ | \ h_{ij}$ when $i
\geq j$.  Forming the matrix $H = (h_{ij}) \in R_p$, we have
$\psi(H) = M$ by construction, and this proves that $\psi$ is
surjective.

Finally, we need to show that $\psi$ is a ring homomorphism.
Clearly, from the definition, $\psi(I) = \text{id}_{E_p}$, and
also $\psi(A+B) = \psi(A)+\psi(B)$. If $A,B \in R_p$, then a
straightforward calculation reveals that $\psi(AB)$ is the
endomorphism composition $\psi(A) \circ \psi(B)$ by the properties
of matrix multiplication. This completes the proof.
\end{proof}

Given this description of End$(H_p)$, one can characterize those
endomorphisms giving rise to elements in Aut$(H_p)$.  Before
beginning this discussion, let us first calculate the kernel of
the map $\psi$ defined in Theorem \ref{endquotientthm}.

\begin{lem}\label{charkerlem}
The kernel of $\psi$ is given by the set of matrices $A =
(a_{ij}) \in R_p$ such that $p^{e_i} \ | \ a_{ij}$ for all $i,j$.
\end{lem}

\begin{proof}
As before, let $w_j = (0,\ldots,g_j,\ldots,0)^T \in H_p$ be the
vector with $g_j$ in the $j$th component and zeroes everywhere
else.  If $A = (a_{ij}) \in R_p$ has the property that each
$a_{ij}$ is divisible by $p^{e_i}$, then
\[\psi(A)w_j = (\pi_1(a_{1j}),\ldots,\pi_n(a_{nj})) = 0.\]  In
particular, since each $h \in H_p$ is a $\mathbb Z$-linear
combination of the $w_j$, it follows that $\psi(A)h = 0$ for all
$h \in H_p$.  This proves that $A \in \ker \psi$.

Conversely, suppose that $A = (a_{ij}) \in \ker \psi$, so that
$\psi(A)w_j = 0$ for each $w_j$.  Then, from the above
calculation, each $a_{ij}$ is divisible by $p^{e_i}$.  This proves
the lemma.
\end{proof}

Theorem \ref{endquotientthm} and Lemma \ref{charkerlem} together
give an explicit characterization of the ring $\text{End}(H_p)$ as
a quotient $R_p/ \ker \psi$.  Following this discussion, we now
calculate the units $\text{Aut}(H_p)$.  The only additional tool
that we require is the following fact from elementary matrix
theory.
\begin{lem}\label{adjugatelem}
Let $A \in \mathbb Z^{n \times n}$ with $\det(A) \neq 0$. Then
there exists a unique matrix $B \in \mathbb Q^{n \times n}$
$($called the \textit{adjugate} of $A$$)$ such that $AB = BA =
\det(A)I$, and moreover $B$ has integer entries.
\end{lem}

Writing $\mathbb F_p$ for the field $\mathbb Z/p\mathbb Z$, the
following is a complete description of $\text{Aut}(H_p)$.

\begin{thm}\label{autdescripthm}
An endomorphism $M = \psi(A)$ is an automorphism if and only if $A
\ (\text{\rm{mod} } p) \in \text{\rm{GL}}_n(\mathbb{F}_p)$.
\end{thm}

\begin{proof}
We begin with a short interlude. Fix a matrix $A \in R_p$ with
$\det(A) \neq 0$. Lemma \ref{adjugatelem} tells us that there
exists a matrix $B \in \mathbb Z^{n \times n}$ such that
$AB = BA = \text{det}(A)I$.  We would like to show that $B$ is
actually an element of $R_p$.  For the proof, express $A =
PA'P^{-1}$ for some $A' \in \mathbb Z^{n \times n}$, and let $B'
\in \mathbb Z^{n \times n}$ be such that $A'B' = B'A' =
\text{det}(A')I$ (again using Lemma \ref{adjugatelem}). Notice
that $\text{det}(A) = \text{det}(A')$.  Let $C = PB'P^{-1}$ 
and observe that \[AC = PA'B'P^{-1} = \text{det}(A)I =
PB'A'P^{-1} = CA.\] By the uniqueness of $B$ from the lemma, it follows
that $B = C = PB'P^{-1}$, and thus $B$ is in $R_p$, as desired.

Returning to the proof of the theorem ($\Leftarrow$), suppose that $p \nmid
\det(A)$ (so that $A \pmod p \in \text{\rm{GL}}_n(\mathbb{F}_p)$),
and let $s \in \mathbb Z$ be such that $s$ is the inverse of
$\det(A)$ modulo $p^{e_n}$ (such an integer $s$ exists since
gcd$(\det(A),p^{e_n}) = 1$).  Notice that we also have $\det(A)
\cdot s \equiv 1 \ (\text{mod } p^{e_j})$ whenever $1 \leq j \leq
n$.  Let $B$ be the adjugate of $A$ as in Lemma \ref{adjugatelem}.
We now define an element of $R_p$, \[A^{(-1)} := s \cdot B,\]
whose image under $\psi$ is the inverse of the endomorphism
represented by $A$:  \[\psi(A^{(-1)}A) = \psi(AA^{(-1)}) = \psi(s
\cdot \det(A)I) = \text{id}_{E_p}.\]  This proves that $\psi(A)
\in \text{Aut}(H_p)$.

Conversely, if $\psi(A) = M$ and $\psi(C) = M^{-1} \in
\text{End}(H_p)$ exists, then \[\psi(AC-I) =
\psi(AC)-\text{id}_{E_p} = 0.\]  Hence, $AC-I \in \ker \psi$. From
the kernel calculation in Lemma \ref{charkerlem}, it follows that
$p \ | \ AC-I$ (entrywise), and so $AC \equiv I \pmod p$.
Therefore, \[1 \equiv \det(AC) \equiv \det(A)\det(C) \pmod p.\] In
particular, $p \nmid \det(A)$, and the theorem follows.
\end{proof}

As a simple application of the above discussion, consider the case
when $e_i = 1$ for $i = 1,\ldots,n$.  Here, $H_p$ can be viewed as
the familiar vector space $\mathbb{F}_p^n$ and $\text{End}(H_p)$
is isomorphic to the ring $\text{M}_n(\mathbb{F}_p)$ of $n \times n$
matrices with coefficients in the field $\mathbb{F}_p$. Theorem
\ref{autdescripthm} is then simply the statement that
Aut($H_p$) corresponds to the set of invertible matrices
$\text{GL}_n(\mathbb{F}_p)$.

\section{Counting the Automorphisms of $H_p$}

To further convince the reader of the usefulness of Theorem
\ref{autdescripthm}, we will briefly explain how to count the number of elements in
$\text{Aut}(H_p)$ using our characterization.  Appealing to Lemma \ref{autproductlemma}, one
then finds an explicit formula for the number of automorphisms of
any finite Abelian group.  The calculation proceeds in two stages: 
$(1)$ finding all elements of $\text{GL}_n(\mathbb{F}_p)$ that can 
be extended to a matrix $A \in R_p$ that represents an endomorphism, 
and then $(2)$ calculating all the distinct ways of extending such an element to
an endomorphism.

Define the following $2n$ numbers: \[d_k = \max \{l : e_l = e_k\},
\ c_k = \min \{l : e_l = e_k\}.\]  Since $e_k = e_k$, we have $d_k
\geq k$ and $c_k \leq k$.  We need to find all $M \in
\text{GL}_n(\mathbb F_p)$ of the form \[ M = \left[
{\begin{array}{*{20}c}
   {m_{11} } & {m_{12} } &  \cdots  & {m_{1n} }  \\
    \vdots  & {} & {} & {}  \\
   {m_{d_1 1} } & {} & {} & {}  \\
   {} & {m_{d_2 2} } & {} & {}  \\
   {} & {} &  \ddots  & {}  \\
   0 & {} & {} & {m_{d_n n} }  \\
 \end{array} } \right] = \left[ {\begin{array}{*{20}c}
   {m_{1c_1 } } & {} & {} & {} & {} & *  \\
   {} & {m_{2c_2 } } & {} & {} & {} & {}  \\
   {} & {} &  \ddots  & {} & {} & {}  \\
   0 & {} & {} & {m_{nc_n } } &  \cdots  & {m_{nn} }  \\
\end{array} } \right].\]
These number \[\prod_{k=1}^{n}(p^{d_k}-p^{k-1}),\] since we only
need linearly independent columns.  Next, to extend each element
$m_{ij}$ from $\overline{m}_{ij} \in \mathbb Z/p\mathbb Z$ to
$\overline{a}_{ij} \in p^{e_i-e_j}\mathbb Z/p^{e_i}\mathbb Z$ such
that \[a_{ij} \equiv m_{ij} \ (\text{mod } p),\] there are
$p^{e_j}$ ways to do this to the necessary zeroes (i.e., when $e_i
> e_j$), since any element of $p^{e_i-e_j}\mathbb Z/p^{e_i}\mathbb
Z$ will do.  Additionally, there are $p^{e_i-1}$ ways at the not
necessarily zero entries ($e_i \leq e_j$), since we may add any element of
$p \mathbb Z/p^{e_i}\mathbb Z$.  This proves the following result.
\begin{thm}
The Abelian group $H_p = \mathbb Z/p^{e_1}\mathbb Z
\times \cdots \times \mathbb Z/p^{e_n}\mathbb Z$ has
\[|\text{\rm Aut}(H_p)| =  \prod_{k=1}^{n}\left(p^{d_k}-p^{k-1}\right) \prod_{j=1}^{n}(p^{e_j})^{n-d_j}
\prod_{i=1}^{n}(p^{e_i-1})^{n-c_i+1}.\]
\end{thm}

\end{document}